\documentclass[letterpaper, 10 pt, conference]{IEEEtran}  

\IEEEoverridecommandlockouts                              
\usepackage{graphics} 
\usepackage{epsfig} 
\usepackage{cite}

\usepackage{amsmath} 
\usepackage{amssymb}  

\title{\LARGE \bf
A Duality Framework for Stochastic Optimal Control of Complex Systems
}

\author{Andreas A. Malikopoulos, {\itshape{Member, IEEE}}
\thanks{This manuscript has been authored by UT-Battelle, LLC, under contract DE-AC05-00OR22725 with the US Department of Energy. The US government retains and the publisher, by accepting the article for publication, acknowledges that the US government retains a nonexclusive, paid-up, irrevocable, worldwide license to publish or reproduce the published form of this manuscript, or allow others to do so, for US government purposes.}
\thanks{This research was supported by the Laboratory Directed Research and Development Program of Oak Ridge National Laboratory, managed by UT- Battelle, LLC, for DOE. This support is gratefully acknowledged.  }
\thanks{A.A. Malikopoulos is with the Energy \& Transportation Science Division, Oak Ridge National Laboratory, Oak Ridge, TN 37831 USA (phone: 865-946-1529; fax: 865-946-1354; e-mail: 
        {\tt\small andreas@ornl.gov)}.}%
}

\begin{document}

\maketitle
\thispagestyle{empty}
\pagestyle{empty}

\begin{abstract}
We address the problem of minimizing the long-run expected average cost of a complex system consisting of interactive subsystems. We formulate a multiobjective optimization problem of the one-stage expected costs of the subsystems and provide a duality framework to prove that the control policy yielding the Pareto optimal solution minimizes the average cost criterion of the system.  We provide the conditions of existence and a geometric interpretation of the solution. For practical situations with constraints consistent to those studied here, our results imply that the Pareto control policy may be of value when we seek to derive online the optimal control policy in complex systems. 
\end{abstract}

\indent

\section{Introduction}
\subsection{Motivation}
Complex systems consist of diverse entities that interact both in space and time. Referring to something as complex implies that it consists of interdependent entities that are connected with each other and can adapt, i.e., they can respond to their local and global environment \cite{Page20007}. Complex systems are encountered in many applications including sustainable transportation, fusion and other alternative energy strategies, and biological systems. For example, the US electricity grid is one of the world's largest complex systems \cite{DOE2011} consisting of a dynamic collection of diverse, interacting components that can adapt. These components are also interdependent and operate under an enormous range of physical, reliability, economic, social, and political constraints that need to be satisfied over time scales ranging from seconds, for closed-loop control, to decades, for transmission siting and construction. Hybrid electric vehicles (HEVs) and plug-in HEVs is another complex system \cite{Malikopoulos2014b} consisting of various interdependent subsystems, e.g., the internal combustion engine, the electric machines (motor and generator), and the energy storage system (battery), that are connected and adapt appropriately to provide the power demanded by the driver. Another example of complex system is the hybrid distributed power generation system \cite{Prasanna2013a} consisting of wind turbines, photovoltaic generation, energy storage, and the relevant energy conversion control.

Stochastic optimal control of complex systems is a ubiquitous task in engineering. The problem is formulated as sequential decision-making under uncertainty where a controller is faced with the task to select control actions in several time steps to achieve long-term goals efficiently. While the nature of these problems may vary widely, their underlying structure is similar and has two principal features: a discrete-time dynamic system whose state evolves according to given transition probabilities that depend on a decision at each time and a cost function that is additive over time. The objective is to derive an optimal policy that minimizes the long-run expected average cost criterion. 

Mathematically, the average cost criterion is prominent as being complex to analyze compared to others; while other classical criteria lead to rational complete solutions, the long-run cost may not \cite{Arapostathis1993}. The average cost criterion in Markov chains with finite state and action spaces is well understood  \cite{Varaiya1978, Berstekas2007, Kushner1971,Kumar1986,Howard1960,Doob1990,Malikopoulos2011}. Dynamic programming (DP) \cite{Bellman1957} has been widely employed as the principal method for analysis of these problems\cite{White1963,Bather1973,Bather1973a,Bather1973b,Hubner1977,Federgruen1978,Bertsekas1998}. A significant amount of work has  focused on inventory problems using linear programming \cite{Manne1960,Wagner1960}, which has been widely used as an alternative to DP method \cite{Derman1962,Derman1965,Hordijk1979,Hordijk1984,Ross1989a,Ross1989b,Lasserre1994,Zadorojniy2006}. Policy iteration \cite{Howard1960} has been another method to address problems considering the average cost criterion  \cite{Miller1969,Chang2007} by adjusting the policy of the system directly rather than using value iteration to derive it. Various other methods proposed in the literature have used matrix decomposition \cite{Lamond1989}, quadratic programming for multiple costs \cite{Ghosh1990}, learning algorithms \cite{Ren2001,Abounadi2001}, decentralized methods \cite{Chang2009}, and the risk-sensitive criterion \cite{Cavazos-Cadena2009}.

Despite the significant progress in optimization and control methods within the last decades current techniques, in some instances, may be computationally impractical for online optimal control of large-scale complex systems \cite{DOE2006}. One possible approach for ameliorating this difficulty is to develop the framework that exploits the structure of the system interconnections and narrow the range of acceptable solutions.

In this paper, we seek to establish a rigorous framework for the analysis and stochastic optimization of complex systems that will permit online implementation of the optimal control policy with respect to the long-run expected average cost criterion. The contributions of this paper are (1) the development of a duality framework for the analysis and stochastic optimization of complex systems that can be used to derive the optimal control policy; (2) the formulation and solution of a multiobjective optimization problem of the one-stage expected costs of all interactive subsystems yielding an equilibrium operating point among the subsystems that minimizes the long-run expected average cost of the system; and (3) the geometric interpretation of the solution and the formation of the conditions under which the optimal control policy exists.

The remainder of the paper proceeds as follows. In Section II, we introduce our notation and formulate the problem. In Section III, we develop a multiobjective optimization framework to address the problem and introduce the Pareto control policy. In Section IV, we show that the Pareto control policy minimizes the long-run expected average cost criterion. Finally, we present illustrative examples in Section V and concluding remarks in Section VI.

\section{System Model and Problem Formulation}
\subsection{Notation}
We denote random variables with upper case letters, and their realization with lower case letters, e.g., for a random variable $X$, $x$ denotes its realization. Subscripts denote time, and subscripts in parentheses denote subsystems; for example, $X_{t(i)}$ denotes the random variable of the subsystem $i$ at time $t$, and $x_{(i)}$ its realization. The shorthand notation $X_{t(1:N)}$ denotes the vector of random variables $\big(X_{t(1)}, X_{t(2)},\cdots,X_{t(N)}\big)$ and $x_{(1:N)}$ denotes the vector of their realization $\big(x_{(1)}, x_{(2)},\cdots,x_{(N)}\big)$. $\mathbb{P}(\cdot)$ is the transition probability matrix, and $\mathbb{E}[\cdot]$ is the corresponding expectation of a random variable. For a control policy $\bf{\pi},$ we use $\mathbb{P}^\pi(\cdot)$, $\mathbb{E}^\pi[\cdot]$ and $\beta^{\pi}$  to denote that the transition probability matrix, expectation and stationary distribution depend on the choice of the control policy $\pi$. 

\subsection{The System Model}
We consider a system consisting of $N$ subsystems. The subsystems interact with each other and their environment. At time $t, t=1,2,\cdots,T$, the state of each subsystem $i, X_{t(i)}$, takes values in a finite state space $\mathcal{S}_{(i)}$, which is a metric space. For each subsystem $i$, we also consider a finite control space $\mathcal{U}_{(i)}$, which is also a metric space,  from which control actions, $U_{t(i)}$, are chosen. 

The initial state of the system $X_{0(1:N)}$ is a random variable taking values in the system's state space, $\mathcal{S}=\prod_{i=1}^N \mathcal{S}_{(i)}$. The evolution of the state is imposed by the discrete-time equation

\begin{equation}
X_{t+1(1:N)}=f(X_{t(1:N)},U_{t(1:N)},W_{t(1:N)}),
\label{eq:1}
\end{equation}
where $W_{t(1:N)}$ is the input from the environment. The system state can be completely observed. 

In our formulation, a state-dependent constraint is incorporated; that is, for each realization of the state of the subsystem $i$, $X_{t(i)}=x_{(i)}$, there is a nonempty and closed set $\mathcal{C}(x_{(i)}):=\big\{u_{(i)}| X_{t(i)}=x_{(i)} \big\} \subset\mathcal{U}_{(i)}$ of feasible control actions when the system is in state $x_{(i)}$. For each subsystem $i$, we denote the set of admissible state/action pairs
\begin{equation}
\Gamma_{(i)} \colon =\{(x_{(i)},u_{(i)} )\vert x_{(i)}\in\mathcal{S}_{(i)}  \text{ and  } u_{(i)} \in \mathcal{C}(x_{(i)}) \}.
\label{eq:4}
\end{equation}
The set of admissible state/action pairs for the system is 
$$
\Gamma \colon=\prod_{i=1}^N \Gamma_{(i)}  =\{(x_{(1:N)},u_{(1:N)}) )\vert x_{(1:N)}\in\mathcal{S}  
$$
\begin{equation}
\text{ and  } u_{(1:N)} \in \mathcal{C}(x_{(1:N)}) \},
\label{eq:4b}
\end{equation}
where $\mathcal{C}(x_{(1:N)}) =\prod_{i=1}^N \mathcal{C}_{(i)}(x_{(i)})$.

For each state of the system $X_{t(1:N)}=x_{(1:N)}$, we define the functions $\mu: \mathcal{S} \to  \mathcal{U}$,  where  $\mathcal{U}=\prod_{i=1}^N \mathcal{U}_{(i)}$, that map the state space to the control action space defined as the control law. When the system is at state $X_{t(1:N)}=x_{(1:N)}$, the controller chooses action according to the control law $u_{(1:N)} = \mu\big(x_{(1:N)}\big)$. 

\textit{Definition 1:} Each sequence of the functions $\mu$ is defined as a stationary control policy of the system
\begin{equation}
\pi \colon =\big(\mu(1),\mu(2),\cdots,\mu(|\mathcal{S}|)\big),
\label{eq:5}
\end{equation}
where $|\mathcal{S}|$ is the cardinality of the system's state space $\mathcal{S}$. 

Let $\Pi$ denote the set of the collection of the stationary control policies 
\begin{equation}
\Pi \colon =\bigg\{\pi \vert \pi= \big(\mu(1),\mu(2),\cdots,\mu(|\mathcal{S}|) \bigg\}.
\label{eq:6}
\end{equation}

The stationary control policy $\pi$ operates as follows. Associated with each state $X_{t(1:N)}=x_{(1:N)}$ is the function $\mu\big(x_{(1:N)}\big)\in\mathcal{C}(x_{(1:N)})$. If at any time the controller finds the system in state $x_{(1:N)}$, then the controller always chooses the action based on the function $\mu\big(x_{(1:N)}\big)$. A stationary policy depends on the history of the process only through the current state, and thus to implement it, the controller only needs to know the current state of the system. The advantages for implementation of a stationary policy are apparent as it requires the storage of less information than required to implement a general policy. 

At each stage $t$, the controller observes the state of the system,  $X_{t(1:N)}=x_{(1:N)}\in\mathcal{S}$, and an action, $u_{t(1:N)}=\mu(X_{t(1:N)})$, is realized from the feasible set of actions at that state. At the same stage $t$, an uncertainty, $W_{t(1:N)}$, is incorporated in the system. At the next stage, $t+1$, the system transits to the state $X_{t+1(1:N)} =x_{(1:N)}'\in\mathcal{S}$ and a transition cost for each subsystem $i$ , $c_{t(i)}\big(X_{t+1(i)} | X_{t(i)},U_{t(i)}\big)$, where $c_{t(i)} :\mathcal{S}_{(i)}\times \mathcal{C}(x_{(i)}) \times \mathcal{S}_{(i)} \to \mathbb{R}$, and for the system, $c_{t}\big(X_{t+1(1:N)} | X_{t(1:N)},U_{t(1:N)}\big)$, where $c_{t}: \mathcal{S}\times \mathcal{C}(x_{(1:N)}) \times \mathcal{S} \to \mathbb{R}$, are incurred.

\subsection{Assumptions}
In the model described above, we consider the following assumptions:

(A1) There exists $\mu$ such that the graph of $\mu$ is included in $\Gamma$.

(A2)  The input from the uncertainty $W_{t(1:N)}$ is a sequence of independent random variables, independent of the initial state $X_{0(1:N)}$, and takes values in the finite set $\mathcal{W}$.

(A3) For each stationary control policy $\pi$, the Markov chain $\left\{ X_{t(1:N)} |t=1,2,\cdots\right\}$ has a unique probability distribution (row vector). 

(A4) The one-stage expected cost of the system, $k_{t}^{\pi} \colon \Gamma \to \mathbb{R},$ 
$$k_{t}^{\pi}\big(X_{t(1:N)},U_{t(1:N)}\big)=$$
$$\sum_{x_{(1:N)}'\in\mathcal{S}} P(X_{t+1(1:N)}=x_{(1:N)}'|X_{t(1:N)}=x_{(1:N)},U_{t(1:N)})\cdot 
$$
$$c_{t}\big(X_{{t+1}(1:N)}=x_{(1:N)}' | X_{t(1:N)}=x_{(1:N)},U_{t(1:N)}\big),$$
is a continuous function of the one-stage costs of the subsystems and it is uniformly bounded.

(A5) The control action realized at each subsystem doesn't affect the transition probability matrix of the other subsystems. 

We briefly comment on the above assumptions. A1 ensures that the set of the collection of the stationary control policies, $\Pi$, is nonempty. A2 imposes a condition yielding that the state $X_{t+1(1:N)}$ depends only on $X_{t(1:N)}$ and $U_{t(1:N)}$. Namely, the evolution of the state is a Markov chain \cite{Kumar1986}. A3 implies that for each stationary policy $\pi\in\Pi$,  there is a unique probability distribution (row vector) 
$\beta^\pi=\big(\beta(1)^{\pi},\beta(2)^{\pi},\cdots,\beta(k)^{\pi},\cdots,\beta(|\mathcal{S}|)^{\pi}\big)$, with 
$\sum_{k=1}^{|\mathcal{S}|}\beta(k)^\pi=1$ \cite[p. 227]{Grimmett2001} such that
$\beta^\pi=\beta^\pi\cdot \mathbb{P}^\pi$. Under this assumption, it is known \cite[p. 175]{Ross1995} that
\begin{equation}
\lim_{T\to\infty}\frac{1}{T+1}\sum_{t=0}^{T}\left[\mathbb{P}^\pi \right]^t=\mathbf{1}\cdot\beta^\pi,
\label{eq:lim}
\end{equation}
where $\mathbb{P^\pi}$ is the transition probability matrix and $\mathbf{1}=[1,1,\cdots,1]^T$. A4 imposes that the interaction of the subsystems has an impact on one-stage expected cost of the system. Finally, A5 implies that the subsystems evolve independently.

\subsection{Problem Formulation}
We are concerned with deriving a stationary optimal control policy $\pi$ to minimize the long-run expected average cost of the system
\begin{equation}
J(\pi)=\lim_{T\to\infty}\frac{1}{T+1}\mathbb{E}^\pi \left [\sum_{0}^{T} k_{t}^{\pi}\big(X_{t(1:N)},U_{t(1:N)}\big)\right].
\label{eq:8a}
\end{equation}
Since for each control policy the Markov chain has a unique probability distribution (A3), it follows that the limit in \eqref{eq:8a} exists. Substituting \eqref{eq:lim} into \eqref{eq:8a} shows that the long-run average cost, $J(\pi)$, does not depend on the initial state $X_{0(1:N)}$ and is given simply as
\begin{equation}
J(\pi)=\beta^\pi\cdot k^\pi,
\label{eq:10}
\end{equation}
where $\beta^\pi$ the the stationary probability distribution of the entire system and
\begin{equation}
k^\pi=\bigg(k^\pi_{t}\big(1,U_{t(1:N)}\big),k^\pi_{t}\big(2,U_{t(1:N)}\big),\cdots k^\pi_{t}\big(|\mathcal{S}|,U_{t(1:N)}\big)\bigg)^T,  
\label{eq:11}
\end{equation}
is the column vector of the system's one-stage expected cost. 

Various methods that discussed in the Introduction can be used to solve \eqref{eq:8a} or \eqref{eq:10} offline and derive the optimal control policy that minimizes the long-run expected average cost $J$ of the system. In this paper, we seek the theoretical framework that will implement the optimal control policy online while the subsystems interact with each other.  The intention here is to identify an equilibrium operating point among the subsystems; if the systems operate at this equilibrium, then the average cost of the system will be minimized.

\section{Multiobjective Optimization Analysis}
\subsection{Pareto Control Policy}
To identify an equilibrium operating point among the subsystems we formulate a multiobjective optimization problem for the one-stage cost of the subsystems. Let's consider the function $f\colon \Gamma \to \mathbb{R}^N$, 
$$
f =\bigg(k^\pi_{t(1)}\big(X_{t(1:N)},U_{t(1:N)}\big), k^\pi_{t(2)}\big(X_{t(1:N)},U_{t(1:N)}\big),...,
$$
\begin{equation}
k^\pi_{t(N)}\big(X_{t(1:N)},U_{t(1:N)}\big)\bigg), 
\end{equation}
where $k^\pi_{t(i)}\big(X_{t(1:N)},U_{t(1:N)}\big)$ is the one-stage expected cost for each subsystem $i$ and the following multiobjective optimization problem
$$
\min_{U_{t(1:N)}\in\mathcal{C}(x_{(1:N)})} \bigg(k^\pi_{t(1)}\big(X_{t(1:N)},U_{t(1:N)}\big), 
$$
\begin{equation}
k^\pi_{t(2)}\big(X_{t(1:N)},U_{t(1:N)}\big),...,k^\pi_{t(N)}\big(X_{t(1:N)},U_{t(1:N)}\big)\bigg). \label{eq:pareto4a}
\end{equation}
The result of the problem \eqref{eq:pareto4a} is called Pareto efficiency. In a Pareto efficiency allocation among agents, no one can be made better without making at least one other agent worse. The following result provides the conditions that the Pareto efficiency exists. 

\textit{Proposition 1} \cite{Ehjrgott2005}: Let $\Gamma$ be a nonempty and compact set, and the one-stage expected cost for each subsystem $i$, $k_{t(i)}^{\pi}(X_{t(i)}, U_{t(i)})\colon\Gamma\to\mathbb{R}$, be lower semicontinuous for all $i=1,\cdots,N$. Then the Pareto efficiency is not empty.

In our problem, the set of admissible state/action pairs, $\Gamma$, is a nonempty compact set (A1). Furthermore, the one-stage expected cost for each subsystem $i$, $k_{t(i)}^{\pi}(X_{t(i)}, U_{t(i)})$, is a continuous function (A4). Consequently, the Pareto efficiency exists. 

\textit{Definition 2:}
The Pareto control policy $\pi^o$ is defined as the policy that yields the minimum one-stage expected cost of the system, $k^{\pi^o}_{t}(X_{t(1:N)}, U^o_{t(1:N)})$, at each realization of the system state $X_{t(1:N)}=x_{(1:N)}$.

\subsection{Impact of the Pareto Control Policy on the System's Expected Cost}
To simplify notation, in the rest of the paper the one-stage expected cost of each subsystem $i$, $k^\pi_{t(i)}\big(X_{t(1:N)},U_{t(1:N)}\big)$, and the one-stage expected cost of the system, $k^\pi_{t}\big(X_{t(1:N)},U_{t(1:N)}\big)$, incurred when the system operates under the control policy $\pi$, will be denoted by $k^\pi_{t(i)}$ and $k^\pi_{t}$ respectively.

\textit{Definition 3:}
In a system consisting of $N$ interactive subsystems, the group of subsystems whose expected costs are a decreasing function with respect to the cost of the system is defined as the \textit{minor} group.

\textit{Definition 4:} 
In a system consisting of $N$ interactive subsystems, the group of subsystems whose expected costs are an increasing function with respect to the cost of the system is defined as the \textit{principal} group.

Without loss of generality, we assume that the minor group consists of the subsystems $1,2,\cdots,m, m\in\mathbb{N}$, and the principal group consists of the subsystems $m+1,\cdots,N$. 
Thus, since the one-stage expected cost of the system is a function $\delta$ of the one-stage cost of the subsystems (A4),
\begin{equation}\label{eq:k}
k^{\pi}_{t}=\delta\big(k^{\pi}_{t(1)},k^{\pi}_{t(2)},\cdots,k^{\pi}_{t(N)} \big),
\end{equation}
from Definition 3, for each subsystem $i$ in the minor group and for any two control policies $\pi,\pi' \in\Pi$ such that $k^{\pi}_{t(i)}\le k^{\pi'}_{t(i)},$ if we fix the one-stage cost of the other subsystems in both minor and principal groups we have
\begin{equation}\label{eq:group1}
k^{\pi}_{t}=\delta\big(\cdots,k^{\pi}_{t(i)},\cdots\big) \ge k^{\pi'}_{t}=\delta\big(\cdots,k^{\pi'}_{t(i)},\cdots\big).
\end{equation}
Similarly, from Definition 4, for each subsystem $j$ in the principal group and for any two control policies $\pi,\pi' \in\Pi$ such that $k^{\pi}_{t(j)}\le k^{\pi'}_{t(j)},$ if we fix the one-stage cost of the other subsystems in both minor and principal groups we have
\begin{equation}\label{eq:group2}
k^{\pi}_{t}=\delta\big(\cdots,k^{\pi}_{t(j)},\cdots,\big)\le k^{\pi'}_{t}=\delta\big(\cdots,k^{\pi'}_{t(j)},\cdots, \big).
\end{equation}

\subsubsection{Problem 1}
We consider the special case where the system consists of $N$ subsystems of a minor group only. 

\textit{Proposition 2:}\label{thm:paretothe1} The solution of the following multiobjective optimization problem at each realization of the state $X_{t(1:N)}=x_{(1:N)}$ yields the Pareto control policy of the system.
\begin{equation}
\begin{split}
\max_{U_{t(1:N)}\in\mathcal{C}(x_{(1:N)})} \big(k^\pi_{t(1)}, \cdots, k^\pi_{t(N)} \big) \label{eq:theorem1a}\\
\text{subject to}  \ X_{t(1:N)} \in\mathcal{S}.
\end{split}
\end{equation}

\begin{proof}
Let $u^*_{(1:N)}\in\mathcal{C}(x_{(1:N)})$ be the solution of \eqref{eq:theorem1a} at each realization of the state $X_{t(1:N)}=x_{(1:N)}$ under the control policy $\pi$. Thus, if we operate the system under $\pi$, then $k_{t(i)}^\pi \ge k_{t(i)}^{\pi'}$, $i=1,\dots,N$, for all $\pi'\in\Pi$. Therefore from Definition 3 we have $k_{t}^\pi \le k_{t}^{\pi'}$, for all $\pi'\in\Pi$, and hence from Definition 2 $\pi$ is the Pareto control policy.
\end{proof}


\subsubsection{Problem 2}
We consider the special case where the system consists of $N$ subsystems of a principal group only. 

\textit{Proposition 3:}\label{thm:paretothe1} The solution of the following multiobjective optimization problem at each realization of the state $X_{t(1:N)}=x_{(1:N)}$ yields the Pareto control policy of the system.
\begin{equation}
\begin{split}
\min_{U_{t(1:N)}\in\mathcal{C}(x_{(1:N)})} \big(k^\pi_{t(1)}, \cdots, k^\pi_{t(N)} \big) \label{eq:theorem2a},\\
\text{subject to}  \ X_{t(1:N)} \in\mathcal{S}.
\end{split}
\end{equation}

\begin{proof}
Let $u^*_{(1:N)}\in\mathcal{C}(x_{(1:N)})$ be the solution of \eqref{eq:theorem2a} at each realization of the state $X_{t(1:N)}=x_{(1:N)}$ under the control policy $\pi$.  Thus, if we operate the system under $\pi$, then $k_{t(i)}^\pi \le k_{t(i)}^{\pi'}$, $i=1,\dots,N$, for all $\pi'\in\Pi$. Therefore from Definition 4 we have $k_{t}^\pi \le k_{t}^{\pi'}$, for all $\pi'\in\Pi$, and hence from Definition 2 $\pi$ is the Pareto control policy.
\end{proof}

\subsubsection{Problem 3}
We consider the general case where the system consists of $N$ subsystems of both a minor and principal group.

In this case, to derive the Pareto control policy, we formulate the following  optimization problem for the one-stage cost of the system
\begin{gather}
\min_{U_{t(1:N)}\in\mathcal{C}(x_{(1:N)})} k^\pi_{t} \label{eq:problem3}  \nonumber \\ 
\min_{U_{t(1:N)}\in\mathcal{C}(x_{(1:N)})} \delta\big(k^{\pi}_{t(1)},k^{\pi}_{t(2)},\cdots,k^{\pi}_{t(N)} \big),\label{eq:problem3}\\ 
\text{subject to}  \ X_{t(1:N)} \in\mathcal{S}. \nonumber
\end{gather}
The Pareto control policy is derived by computing at each realization of the system state $X_{t(1:N)}=x_{(1:N)} \in\mathcal{S}$,  the control action $u_{(1:N)}^o$  that yields the minimum one-stage expected cost of the system in \eqref{eq:problem3}.

\section{Duality Framework}

\subsection{Geometric Framework for Duality Analysis}

We  use a geometric framework from duality analysis, referred to as \textit{min common/max crossing point} problems \big(see \cite{Bertsekas2003}, p. 120\big), to show that the Pareto control policy is an optimal control policy that minimizes the long-run expected average cost of the system, and provide a geometric interpretation of the solution. 

\begin{figure}[ht]
  \centering
    \includegraphics[width=3.5 in]{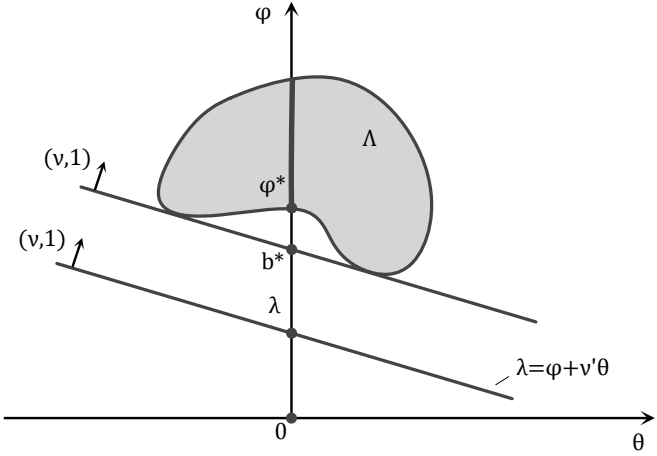} 
      \caption{Geometric framework for duality analysis.}
      \label{fig:2}
\end{figure}

The \textit{min common/max crossing point} framework captures the most essential elements of duality by considering two geometric problems. Let's consider a nonempty subset $\Lambda$ of $\mathbb{R}^{n+1}$ as shown in Fig. \eqref{fig:2}. The axis $\theta$ corresponds to $\mathbb{R}^{n}$ and the axis $\varphi$ corresponds to $\mathbb{R}$. 

The first geometric problem, the \textit{min common point}, seeks to find the minimum value $\varphi^*$ of the subset $\Lambda$ in $\varphi$ axis. The second geometric problem, the \textit{max crossing point}, seeks to find the nonvertical hyperplane that contains $\Lambda$ in its corresponding upper closed half space and crosses $\varphi$ axis at a maximum point $b^*$.

Mathematically, the \textit{min common point} problem can be written as
\begin{equation}\label{eq:pareto10}
\min {\varphi}
\end{equation}
$$\text{subject to}: (0,\varphi)\in\Lambda$$

A nonvertical hyperplane in $\mathbb{R}^{n+1}$ is specified by its normal $(\nu,1)\in\mathbb{R}^{n+1}$, where $\nu\in\mathbb{R}^n$, and a scalar $\lambda\in\mathbb{R}$ as
\begin{equation}\label{eq:pareto11}
\varphi+\nu'\theta=\lambda.
\end{equation}
Such a hyperplane crosses the $(n+1)st$ axis, $\varphi$, at $(0,\lambda)$. The hyperplane contains $\Lambda$ in its upper closed half plane if and only if for all $(\theta, \varphi)\in\Lambda$
\begin{equation}\label{eq:pareto11}
\varphi+\nu'\theta\ge\lambda.
\end{equation}
Similarly
\begin{equation}\label{eq:pareto12}
\inf_{(\theta, \varphi)\in\Lambda}{\big\{\varphi+\nu'\theta\big\}}\ge\lambda.
\end{equation}
Thus the \textit{max crossing point} problem can be written 
\begin{equation}\label{eq:pareto13}
\max \inf_{(\theta, \varphi)\in\Lambda} {\big\{\varphi+\nu'\theta\big\}}
\end{equation}
$$\text{subject to}: \nu\in\mathbb{R}^n$$
The function $b(\nu)=\inf{\big\{\varphi+\nu'\theta\big\}}$ is the dual function.


\textit{Definition 5:} \label{def:IV1} 
If $(\bar{\theta},\bar{\varphi})$ belongs to the closure of $\Lambda$ and for all $(\theta,\varphi)\in\Lambda$,  $\bar{\theta}+\nu'\cdot\bar{\varphi}\le\theta+\nu'\cdot\varphi$, we say that the hyperplane supports $\Lambda$ at $(\bar{\theta},\bar{\varphi})$.

\textit{Proposition 4:} (see \cite{Bertsekas2003}, p. 123\big)
The \textit{max crossing point} of the dual function is less than or equal to the \textit{min common point}, namely $b^*\le\varphi^*$.

\begin{proof}
For all $(\theta,\varphi)\in\Lambda$ and $\nu\in\mathbb{R}^n$ we have
\begin{equation}\label{eq:prop5a}
b(\nu)= \inf_{(\theta, \varphi)\in\Lambda}{\big\{\varphi+\nu'\theta\big\}} \le \inf_{(0, \varphi)\in\Lambda}{\big\{\varphi \big\}}.
\end{equation}
Taking the supremum over $\nu\in\mathbb{R}^n$, we have
\begin{equation}\label{eq:prop5b}
b^*= \sup_{\nu\in\mathbb{R}^n}\inf_{(\theta, \varphi)\in\Lambda}{\big\{\varphi+\nu'\theta\big\}} \le \varphi^*=\sup_{\nu\in\mathbb{R}^n}\inf_{(0, \varphi)\in\Lambda}{\big\{\varphi \big\}}.
\end{equation}
\end{proof}

\subsection{Strong Duality of the Pareto Control Policy}
We want to investigate the impact of the Pareto control policy on the long-run expected average cost of the system. This will involve characterizing the solution of the Pareto control policy within a duality framework. We recall that $k^\pi$ is the column vector of the system's one-stage expected cost for each state, $1,2,\cdots,|\mathcal{S}|$, under the control policy $\pi=\big(\mu(1),\mu(2),\cdots,\mu(|\mathcal{S}|)\big)$, namely
$$k^\pi=\bigg(k^\pi\big(1,\mu(1)\big),k^\pi\big(2,\mu(2)\big),\cdots,k^\pi\big(|\mathcal{S}|,\mu(|\mathcal{S}|)\big)\bigg)^T. 
$$

We formulate the following problem: 

\begin{equation}\label{eq:pareto14}
\min_{\pi\in\Pi} \lVert k^\pi + \mathbb{M}^\pi\cdot q \rVert
\end{equation}
$$\text{subject to}: \beta^\pi\cdot\mathbb{M}^\pi=0,$$
where $\mathbb{M}^{\pi}= \mathbb{P}^\pi - \text{I}$, $q\in\mathbb{R^{|\mathcal{S}|}}$ such that $\mathbb{M}^\pi \cdot q > 0$, and $\beta^\pi=\big(\beta(1)^{\pi},\beta(2)^{\pi},\cdots,\beta(|\mathcal{S}|)^{\pi}\big)$ is the probability distribution corresponding to the control policy $\pi$.

We refer to this problem as the primal problem, and we denote by $\lVert k^\pi+ \mathbb{M}^\pi\cdot q \rVert^*$ its optimal value.
The Lagrangian function of the above minimization problem is
\begin{equation}\label{eq:pareto16}
L(\pi,\nu)= \lVert k^\pi+ \mathbb{M}^\pi\cdot q \rVert + \big(\beta^\pi\cdot\mathbb{M}^\pi\big)\cdot\nu,
\end{equation}
where $\nu\in\mathbb{R^{|\mathcal{S}|}}$ is the vector of the Lagrange multipliers.

We use the \textit{min common/max crossing point} framework described above to visualize the duality in \eqref{eq:pareto16}. We consider the following set
\begin{equation}\label{eq:pareto17}
\Lambda:=\big\{\big(\beta^\pi\cdot\mathbb{M}^\pi, \lVert k^\pi+ \mathbb{M}^\pi\cdot q \rVert |\pi\in\Pi\big\}.
\end{equation}


\textit{Lemma 1:} The hyperplane with norm $(\nu,1)$ that passes through the vector $\big(\beta^\pi\cdot\mathbb{M}^\pi, \lVert k^\pi + \mathbb{M}^\pi\cdot q\rVert\big)$ intercepts the vertical axis $\varphi$ at the value of $L(\pi,\nu)$.

\begin{proof}
The hyperplane with norm $(\nu,1)$ that passes through $\big(\beta^\pi\cdot\mathbb{M}^\pi, \lVert k^\pi + \mathbb{M}^\pi\cdot q \rVert\big)$ satisfies
\begin{equation}\label{eq:pareto18}
\varphi+\theta' \cdot\nu=\lVert k^\pi+ \mathbb{M}^\pi\cdot q \rVert+\big(\beta^\pi\cdot\mathbb{M}^\pi\big)\cdot\nu=L(\pi,\nu).
\end{equation}
\end{proof}


\textit{Lemma 2:} The hyperplane that passes through $\lVert k^\pi + \mathbb{M}^\pi\cdot q \rVert^*$ supports $\Lambda$.

\begin{proof}
From Lemma 1 we have
\begin{equation}\label{eq:pareto19}
\varphi+\theta' \cdot\nu= \lVert k^\pi + \mathbb{M}^\pi\cdot q\rVert+\big(\beta^\pi\cdot\mathbb{M}^\pi\big)\cdot\nu=\lambda.
\end{equation}
Since for each stationary control policy we have a unique probability distribution (A3),
\begin{equation}\label{eq:pareto20}
\begin{split}
\beta^\pi=\beta^\pi\cdot\mathbb{P}^\pi\Rightarrow
\beta^\pi\cdot\big(\mathbb{P}^\pi-\text{I}\big)=0\Rightarrow\\
\big(\beta^\pi\cdot\mathbb{M}^\pi\big)=0.
\end{split}
\end{equation}
Thus for each control policy $\pi\in\Pi$, the elements of the set $\Lambda$ are located only on the axis $\varphi$, and 
\begin{equation}\label{eq:pareto21}
\lVert k^\pi + \mathbb{M}^\pi\cdot q\rVert =\lambda.
\end{equation}
Thus
\begin{equation}\label{eq:pareto22}
\lVert k^{\pi^*} + \mathbb{M}^{\pi^*}\cdot q\rVert^*=\varphi^* \le \lVert k^\pi + \mathbb{M}^\pi\cdot q \rVert =\lambda, \forall\pi\in\Pi.
\end{equation}
\end{proof}


\textit{Theorem 1:} The Pareto control policy $\pi^o$ is the optimal control policy that minimizes the long-run expected average cost criterion of the system, under the assumption (A3) and (A4).

\begin{proof}
Let
\begin{gather}
\mathbf{1}\cdot \psi= k^{\pi} + \mathbb{M}^{\pi}\cdot q,  \qquad \forall \pi\in\Pi, \label{eq:pareto25}
\end{gather} 
where $\mathbf{1}=\big(1,1,...,1\big)^T$, and $\psi^{\pi}\in\mathbb{R}$. Recall that $q\in\mathbb{R^{|\mathcal{S}|}}$ such that $\mathbb{M}^\pi \cdot q > 0$.

Multiplying the above equation by $\beta^{\pi}=\big(\beta(1)^{\pi},\beta(2)^{\pi},\cdots,\beta(k)^{\pi},\cdots,\beta(|\mathcal{S}|)^{\pi}\big)$ from the left we have
\begin{gather}\label{eq:pareto26}
 \psi^{\pi}= \beta^{\pi}\cdot k^{\pi}+ \beta^{\pi}\cdot\mathbb{M}^{\pi}\cdot q\\
=\beta^{\pi}\cdot k^{\pi} + \beta^{\pi}\cdot \big(\mathbb{P}^{\pi}-\text{I}\big)\cdot q \\
=\beta^{\pi}\cdot k^{\pi} + \beta^{\pi}\cdot \mathbb{P}^{\pi} \cdot q -\beta^{\pi^o}\cdot q\\
=\beta^{\pi}\cdot k^{\pi} + \beta^{\pi}\cdot q -\beta^{\pi}\cdot q=\beta^{\pi}\cdot k^{\pi}
\end{gather}
since $\mathbb{M}^{\pi}= \mathbb{P}^{\pi} - \mathbb{I}$ and $\beta^{\pi}=\beta^{\pi}\cdot \mathbb{P}^{\pi}$.
So from \eqref{eq:10}, $\psi^{\pi}$ is the long-run expected average cost corresponding to the control policy $\pi$.

From the Definition 2 of the Pareto control policy
\begin{equation}\label{eq:pareto26b}
\begin{split}
k^{\pi^o} \le k^\pi, \qquad \forall \pi\in\Pi,
\end{split}
\end{equation}
and since $\mathbb{M}^\pi \cdot q > 0$, \eqref{eq:pareto26b} through \eqref{eq:pareto25} can be written
\begin{equation}\label{eq:pareto29a}
k^{\pi^o} \le k^\pi + \mathbb{M}^\pi \cdot q = \mathbf{1}\cdot\psi^\pi,
\end{equation}
where $\psi^\pi$ is the long-run expected average cost corresponding to any control policy $\pi\in\Pi$.
Multiplying \eqref{eq:pareto29a} by $\beta^{\pi^o}$ from the left we have
\begin{equation}\label{eq:pareto29b}
\psi^{\pi^o}=\beta^{\pi^o} \cdot k^{\pi^o} \le \psi^\pi, \qquad \forall \pi\in\Pi.
\end{equation}
Thus the Pareto control policy is the optimal control policy that minimizes the long-run expected average cost.
\end{proof}

\textit{Theorem 2:} The Pareto control policy $\pi^o$ supports $\Lambda$.
\begin{proof}
From Theorem 1 we have
\begin{equation}
k^{\pi^o} + \mathbb{M}^{\pi^o} \cdot q \le k^\pi + \mathbb{M}^\pi \cdot q, \qquad \forall \pi\in\Pi.
\end{equation}

Hence
\begin{equation} \label{eq:pareto22b}
\lVert k^{\pi^o} + \mathbb{M}^{\pi^o} \cdot q \rVert \le \lVert k^\pi + \mathbb{M}^\pi \cdot q \rVert.
\end{equation}
and from Lemma 2
\begin{equation} \label{eq:pareto22c}
\lVert k^{\pi^*} + \mathbb{M}^{\pi^*}\cdot q\rVert^*=\lVert k^{\pi^o} + \mathbb{M}^{\pi^o} \cdot q \rVert.
\end{equation}
\end{proof}

\textit{Corollary 1:} There is no duality gap in \eqref{eq:pareto16}, and  thus the Pareto control policy $\pi^o$ yields the global optimal solution.

\section{Illustrative Examples}
\subsection{Preliminary Results}
In this section, we provide some results that we need to use for the illustrative examples in the next subsection. We begin by recalling the Kronecker product and its properties (see \cite{Searle2006}, \cite{Horn1994}).

\textit{Definition 6:}
If A is an m-by-n matrix and B is a p-by-q matrix, then the Kronecker product $A \otimes B$ is the mp-by-np block matrix
$$
A \otimes B=\left[\begin{array}{ccc} a_{11}B & \cdots & a_{1n}B \\
                           \vdots & \ddots & \vdots \\
                           a_{m1}B & \cdots & a_{mn}B \end{array} \right].
$$
The next proposition provides an expression of the transition probability of the entire system as a Kronecker product of the transition probabilities of each subsystem.

\textit{Proposition 5\cite{Malikopoulos2015} :} 
Consider $N$ evolving subsystems with corresponding transition probability matrices $\mathbb{P}_{(i)}$, $i=1, \cdots, N$ defined by $\mathbb{P}_{(i)}(X_{t+1(i)}=x'_{(i)}|X_{t(i)}=x_{(i)}, U_{t(i)}=u_{(i)})$. Now consider that the system operates under the control policy $\pi$. Then the transition probability matrix of the entire system satisfies
\begin{equation}\label{eq:trans}
\mathbb{P}^\pi= \mathbb{P}_{(1)}^{\pi}\otimes \mathbb{P}_{(2)}^{\pi}\otimes \cdots \otimes \mathbb{P}_{(N)}^{\pi}.
\end{equation}

\textit{Proposition 6 \cite{Malikopoulos2015}:}
Consider a controlled Markov chain with a unique probability distribution for each control policy $\pi$ (A3) for the entire system and another one for each subsystem. Then the stationary probability of the entire system, $\beta^\pi$, can be expressed as the Kronecker product of each stationary probability of each corresponding subsystem $i$, $\beta^{\pi}_{(i)}$,  $i=1,\cdots,N$, i.e.,
\begin{equation}
\beta^\pi=\beta^{\pi}_{(1)}\otimes \beta^{\pi}_{(2)}\otimes \cdots \otimes \beta^{\pi}_{(N)}. \label{eq:beta}
\end{equation}

\subsection{A System with Subsystems of a Minor Group}
\begin{figure}[ht]
  \centering
    \includegraphics[width=2.5 in]{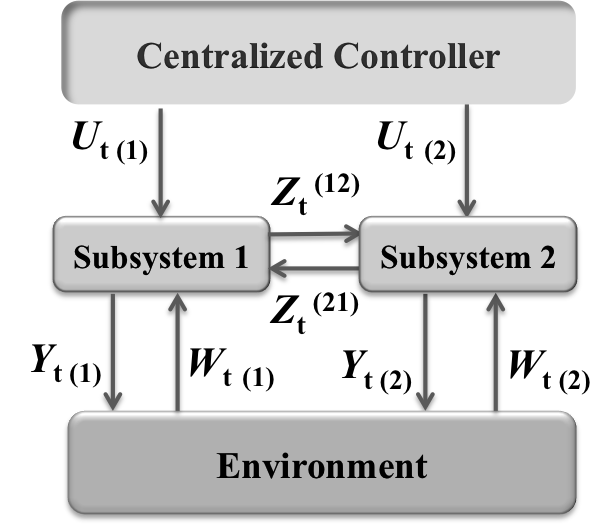} 
      \caption{A System of two subsystems.}
      \label{fig:3}
\end{figure}

We consider a system of two interactive subsystems of a minor group \cite{Malikopoulos2015}, illustrated in Fig. \ref{fig:3}. 

Each subsystem has two states, i.e., $\mathcal{S}_{(i)}=\{1,2\}$, and two control actions $\mathcal{U}_{(i)}=\{a,b\}$. Thus the system has four states $\mathcal{S}=\{1,2, 3, 4\}=\big\{ \left[\begin{array}{cc}1\\1\end{array}\right],\left[\begin{array}{cc}1\\2\end{array}\right],\left[\begin{array}{cc}2\\1\end{array}\right],\left[\begin{array}{cc}2\\2\end{array}\right] \big\}$, and there are sixteen control policies. The transition probability matrices associated with the control policies for the first subsystem are $\mathbb{P}_{(1)}^{\pi^1}=\mathbb{P}_{(1)}^{\pi^2}=\mathbb{P}_{(1)}^{\pi^3}=\mathbb{P}_{(1)}^{\pi^4}=\left[\begin{array}{cc}0.7 & 0.3\\ 0.4 & 0.6\end{array}\right]$, $\mathbb{P}_{(1)}^{\pi^5}=\mathbb{P}_{(1)}^{\pi^6}=\mathbb{P}_{(1)}^{\pi^7}=\mathbb{P}_{(1)}^{\pi^8}=\left[\begin{array}{cc}0.7 & 0.3\\ 0.2 & 0.8\end{array}\right]$, $\mathbb{P}_{(1)}^{\pi^9}=\mathbb{P}_{(1)}^{\pi^{10}}=\mathbb{P}_{(1)}^{\pi^{11}}=\mathbb{P}_{(1)}^{\pi^{12}}=\left[\begin{array}{cc}0.9 & 0.1\\ 0.4 & 0.6\end{array}\right]$,  and $\mathbb{P}_{(1)}^{\pi^{13}}=\mathbb{P}_{(1)}^{\pi^{14}}=\mathbb{P}_{(1)}^{\pi^{15}}=\mathbb{P}_{(1)}^{\pi^{16}}=\left[\begin{array}{cc}0.9 & 0.1\\ 0.2 & 0.8\end{array}\right]$. Similarly, the transition probability matrices for the second subsystem are $\mathbb{P}_{(2)}^{\pi^1}=\mathbb{P}_{(2)}^{\pi^5}=\mathbb{P}_{(2)}^{\pi^9}=\mathbb{P}_{(2)}^{\pi^{13}}=\left[\begin{array}{cc}0.5 & 0.5\\0.45 &0.55\end{array}\right]$, $\mathbb{P}_{(2)}^{\pi^2}=\mathbb{P}_{(2)}^{\pi^6}=\mathbb{P}_{(2)}^{\pi^{10}}=\mathbb{P}_{(2)}^{\pi^{14}}=\left[\begin{array}{cc}0.5 & 0.5\\ 0.3 & 0.7\end{array}\right]$, $\mathbb{P}_{(2)}^{\pi^3}=\mathbb{P}_{(2)}^{\pi^7}=\mathbb{P}_{(2)}^{\pi^{11}}=\mathbb{P}_{(2)}^{\pi^{15}}=\left[\begin{array}{cc}0.6 & 0.4\\ 0.45 & 0.55\end{array}\right]$, and $\mathbb{P}_{(2)}^{\pi^4}=\mathbb{P}_{(2)}^{\pi^8}=\mathbb{P}_{(2)}^{\pi^{12}}=\mathbb{P}_{(2)}^{\pi^{16}}=\left[\begin{array}{cc}0.6 & 0.4\\0.3 &0.7\end{array}\right]$. 

The output for each subsystem with respect to each control policy is given by four $2\times2$ matrices as we have two states and two actions for each subsystem. For the first subsystem corresponding to each control policy the output is given: $Y_{t(1)}^{\pi^1}=Y_{t(1)}^{\pi^2}=Y_{t(1)}^{\pi^3}=Y_{t(1)}^{\pi^4}=\left[\begin{array}{cc}4.8 & 4.0\\ 5.6 & 9.6\end{array}\right]$, $Y_{t(1)}^{\pi^5}=Y_{t(1)}^{\pi^6}=Y_{t(1)}^{\pi^7}=Y_{t(1)}^{\pi^8}=\left[\begin{array}{cc}4.8 & 4.0\\11.2 & 10.4\end{array}\right]$, $Y_{t(1)}^{\pi^9}=Y_{t(1)}^{\pi^{10}}=Y_{t(1)}^{\pi^{11}}=Y_{t(1)}^{\pi^{12}}=\left[\begin{array}{cc}8.0 & 6.4\\5.6 & 9.6\end{array}\right]$, and $Y_{t(1)}^{\pi^{13}}=Y_{t(1)}^{\pi^{14}}=Y_{t(1)}^{\pi^{15}}=Y_{t(1)}^{\pi^{16}}=\left[\begin{array}{cc}8.0 & 6.4\\11.2 & 10.4\end{array}\right]$. The output of the second subsystem with respect to each control policy is $Y_{t(2)}^{\pi^1}=Y_{t(2)}^{\pi^5}=Y_{t(2)}^{\pi^9}=Y_{t(2)}^{\pi^{13}}=\left[\begin{array}{cc}4.9 & 4.2\\ 6.3 & 7.0\end{array}\right]$, $Y_{t(2)}^{\pi^2}=Y_{t(2)}^{\pi^6}=Y_{t(2)}^{\pi^{10}}=Y_{t(2)}^{\pi^{14}}=\left[\begin{array}{cc}4.9 & 4.2\\7.7 & 9.8\end{array}\right]$, $Y_{t(2)}^{\pi^3}=Y_{t(2)}^{\pi^7}=Y_{t(2)}^{\pi^{11}}=Y_{t(2)}^{\pi^{15}}=\left[\begin{array}{cc}6.3 & 8.4\\6.3 & 7.0\end{array}\right]$, and $Y_{t(2)}^{\pi^4}=Y_{t(2)}^{\pi^8}=Y_{t(2)}^{\pi^{12}}=Y_{t(2)}^{\pi^{16}}=\left[\begin{array}{cc}6.3 & 8.4\\7.7 & 9.8\end{array}\right]$.

We assume that $25\%$ of the subsystem's output goes to subsystem $2$, i.e., $Z_{t}^{(12)}=0.25\cdot Y_{t(1)}$ and also $43\%$ percent of the subsystem's output goes to subsystem $1$, i.e., $Z_{t}^{(21)}=0.43\cdot Y_{t(2)}.$ The input for each subsystem is $W_{t(1)}=15$ and $W_{t(2)}=16$ respectively. Furthermore, we assume that the transition cost for each subsystem is given by
\begin{align}
c_{t(1)}\big(X_{t(1)},U_{t(1)}\big)=\frac{W_{t(1)}+Z_{t}^{(21)}}{Y_{t(1)}+Z_{t}^{(12)}}, \label{eq:13}
\end{align} 
and
\begin{align}
c_{t(2)}\big(X_{t(2)},U_{t(2)}\big)=\frac{W_{t(2)}+Z_{t}^{(12)}}{Y_{t(2)}+Z_{t}^{(12)}} \label{eq:14} 
\end{align} 
respectively. The transition cost for the entire system is given by 
\begin{align}
c_{t}\big(X_{t(1:2)},U_{t(1:2)}\big)=\frac{W_{t(1)}+W_{t(2)}}{Y_{t(1)}+Y_{t(2)}}. \label{eq:15}
\end{align}

The transition cost matrix for each subsystem and for the entire system is a $4\times4$ matrix since we have four states in total (two for each subsystem), and the cost depends on each state and control action. For example, if we want to compute the transition cost matrices for each subsystem, $\mathbb{C}^{\pi^1}_{(1)}, \mathbb{C}^{\pi^1}_{(2)}$, and for the system, $\mathbb{C}^{\pi^1}$,  when the system operates under the control policy $\pi^1$, substituting $W_{t(1)},W_{t(2)}, Y_{t(1)}^{\pi^1}, Z_{t}^{(12)}, Y_{t(2)}^{\pi^1}, Z_{t}^{(21)}$ in \eqref{eq:13},  \eqref{eq:14}, and \eqref{eq:15} yields
\begin{equation*}
\mathbb{C}^{\pi^1}_{(1)}=\left[\begin{array}{cccc}2.85 & 2.80 & 3.42 & 3.36\\ 2.95 & 3.00 & 3.54 & 3.60\\ 2.44 & 2.40 & 1.43 & 1.40\\ 2.53 & 2.57 & 1.48 & 1.50 \end{array}\right],
\end{equation*}

\begin{equation*}
\mathbb{C}^{\pi^1}_{(2)}=\left[\begin{array}{cccc}2.45 & 2.87 & 2.43 & 2.83\\ 1.56 &  1.23 & 1.55 & 1.21\\ 2.69 & 3.13 & 2.66 & 3.10\\  1.71 & 1.34 & 1.69 &1.33 \end{array}\right], \text{and}
\end{equation*}

\begin{equation*}
\mathbb{C}^{\pi^1}=\left[\begin{array}{cccc}3.19 & 3.44 & 3.48 & 3.78\\ 2.48 &   2.12 & 2.64 & 2.24\\ 1.92 & 2.01 & 2.02 & 2.12\\  1.64 & 1.47 & 1.71 & 1.53 \end{array}\right].
\end{equation*}
The entry $(3, 2)$ in the transition cost matrices $\mathbb{C}^{\pi^1}_{(1)}, \mathbb{C}^{\pi^1}_{(2)},$ and $\mathbb{C}^{\pi^1}$ corresponds to the costs incurred when the subsystem $1$ resides at state $2$ and transits to state $1$ while the subsystem $2$ resides at state $1$ and transits to state $2$ following the control policy $\pi^1$.

Similar to the cost matrix, the transition probability matrix is also a $4\times 4$ for the four states. When the system operates under the control policy $\pi^1$, the transition probability matrix is given from Proposition 5, i.e., $\mathbb{P}^{\pi^1} = \mathbb{P}_{(1)}^{\pi^1} \otimes \mathbb{P}_{(2)}^{\pi^1}.$ Therefore,
\begin{center}
$\mathbb{P}^{\pi^1}=\left[\begin{array}{cccc} 0.35 & 0.35 & 0.15 & 0.15\\ 0.315 & 0.385 & 0.135 & 0.165\\ 0.2 & 0.2 & 0.3 & 0.3\\ 0.18 & 0.22 & 0.27 & 0.33\end{array}\right]$.
\end{center}

The one-stage expected cost, $k^{\pi}\big(X_{t(1:2)},U_{t(1:2)}\big),$ of each subsystem $i$ is a $4\times1$ vector, and the value of the element $m$ is computed as follows:

$$k_{t(i)}^{\pi}\big(X_{t(1:2)},U_{t(1:2)}\big)$$
\begin{align}
=\sum_{k=1}^{4}[\mathbb{P}^{\pi}]_{mk} \cdot [\mathbb{C}^{\pi}]_{mk}.\label{eq:onecost}
\end{align} 
For example, to compute the one-stage expected cost for subsystem $1$ following the control policy $\pi^1$ we have
$$
k^{\pi^1}_{(1)}\big(X_{t(1:2)},U_{t(1:2)}\big)$$
$$=
\left[\begin{array}{c} \sum_{k=1}^4{[\mathbb{P}^{\pi^1}]_{1k}}[\mathbb{C}^{\pi^1}_{(1)}]_{1k}\\ \sum_{k=1}^4{[\mathbb{P}^{\pi^1}]_{2k}}[\mathbb{C}^{\pi^1}_{(1)}]_{2k}\\ \sum_{k=1}^4{[\mathbb{P}^{\pi}]_{3k}}[\mathbb{C}^{\pi^1}_{(1)}]_{3k}\\ \sum_{k=1}^4{[\mathbb{P}^{\pi^1}]_{4k}}[\mathbb{C}^{\pi^1}_{(1)}]_{4k}\end{array}\right]$$
\begin{align}
=\left[\begin{array}{c}2.9945\\ 3.1562 \\ 1.8170\\ 1.9154\end{array}\right].
\end{align}

The stationary probability distribution is given by \eqref{eq:beta}. For example, the stationary distribution imposed by the control policy $\pi^1$, is $\beta^{\pi^1}=\beta^{\pi^1}_{(1)} \otimes \beta^{\pi^1}_{(2)}=[\begin{array}{cccc}0.2707 & 0.3008 & 0.2030 & 0.2256\end{array}]$. Hence the average cost of subsystem 1 with respect to policy $\pi^1$ is given by \eqref{eq:10}, $J(\pi)=\beta^\pi\cdot k^{\pi^1}_{(1)}=2.5602$. In a similar way we can compute the corresponding one-stage cost vectors and probability distributions for the subsystems $1$, $2$, and the entire system for all $16$ control policies. The average costs for the subsystems and the system corresponding to each control policy are summarized in Tables \ref{ta:de1}, \ref{ta:de2} and \ref{ta:de3}. Each value in the table (reading the table row by row) corresponds to the long-run expected average cost for the control policies from $\pi^1$ to $\pi^{16}$. We note that subsystem 1 reaches its minimum average cost $J_1$ when the policy $\pi^{13}$ is used. For subsystem 2, the optimal cost is attained with the policy $\pi^{4}$. Finally, for the entire system optimality occurs under the control policy $\pi^{16}$ which is the Pareto control policy as it corresponds to the Pareto efficiency one-stage expected cost for each subsystem.

\begin{table}[ht]
\begin{center}
\caption{Long-Run Average Costs for Subsystem 1}
\begin{tabular}{|c|c|c|c|c|}
  \hline
  2.5602 & 2.6712 & 2.6390 & 2.7255 \\
  \hline
  2.0249 & 2.1127 & 2.0872 & 2.1556 \\
  \hline
  1.8029 & 1.8811 & 1.8584 & 1.9193 \\
  \hline
  1.6317 & 1.7025 & 1.6820 & 1.7371 \\
  \hline
\end{tabular}
\label{ta:de1}
\end{center}
\end{table}

\begin{table}[ht]
\begin{center}
\caption{Long-Run Average Costs for Subsystem 2}
\begin{tabular}{|c|c|c|c|c|}
  \hline
2.2511 & 1.8617 & 1.7921 & 1.5235 \\
  \hline
2.3194 & 1.9182 & 1.8464 & 1.5697 \\
  \hline
2.3102 & 1.9106 & 1.8391 & 1.5634 \\
  \hline
2.3383 & 1.9338 & 1.8615 & 1.5825 \\
  \hline
\end{tabular}
\label{ta:de2}
\end{center}
\end{table}

\begin{table}[ht]
\begin{center}
\caption{Long-Run Average Costs for Entire System}
\begin{tabular}{|c|c|c|c|c|}
  \hline
2.7557 & 2.4427 & 2.4607 & 2.2307 \\
  \hline
2.3801 & 2.1328 & 2.1522 & 1.9695 \\
  \hline
2.3178 & 2.0876 & 2.1108 & 1.9398 \\
  \hline
2.1821 & 1.9746 & 1.9977 & 1.8431\\
  \hline
\end{tabular}
\label{ta:de3}
\end{center}
\end{table}

\subsection{A System with Subsystems of a Minor Group with Varying Transition Probability and Cost Matrices}
In this example we use synthetic data to examine the Pareto control policy of the systems. First, we use DP to compute the optimal control policy, denoted by $\pi^*$, that minimizes the average cost of the entire system. We anticipate that the Pareto control policy $\pi^0$ will yield the same result.

Let the subsystems' inputs be $W_{t(1)}=15, W_{t(2)}=16$ as in the previous example. Next, a random output of each subsystem is considered. The total output of the first subsystem 1 associated with action, $a$, is a matrix with random entries distributed according to a uniform distribution, $Y(1,3)$. Similarly, the total output of the same subsystem with respect to action $b$, is a matrix with entries distributed according to $Y(8,10)$. For the second subsystem, the entries of the matrix associated with action $a$ are independent and identically distributed (i.i.d.) $Y(2,4)$, and the ones  associated with action $b$ i.i.d. $Y(9,12)$. 

Next, let $\rho^*\doteq \rho(\pi^*)$ be the map defined as 
\begin{equation}\label{eq:rho}
\rho^{\pi} = \|f-f^s\|,
\end{equation}
where 
\begin{gather}
f= \bigg(k^\pi_{t(1)}\big(X_{t(1:2)},U_{t(1:2)}\big),  k^\pi_{t(2)}\big(X_{t(1:2)},U_{t(1:2)}\big) \bigg),
\end{gather}
and
\begin{gather}
f^s= \bigg(\min_{U_{t(1:2)}\in\mathcal{C}(x_{(1:2)})} k^\pi_{t(1)}\big(X_{t(1:2)},U_{t(1:2)}\big), \nonumber \\
\min_{U_{t(1:2)}\in\mathcal{C}(x_{(1:2)})} k^\pi_{t(2)}\big(X_{t(1:2)},U_{t(1:2)}\big) \bigg).
\end{gather}

\begin{figure}[ht]
  \centering
    \includegraphics[width=3.7 in]{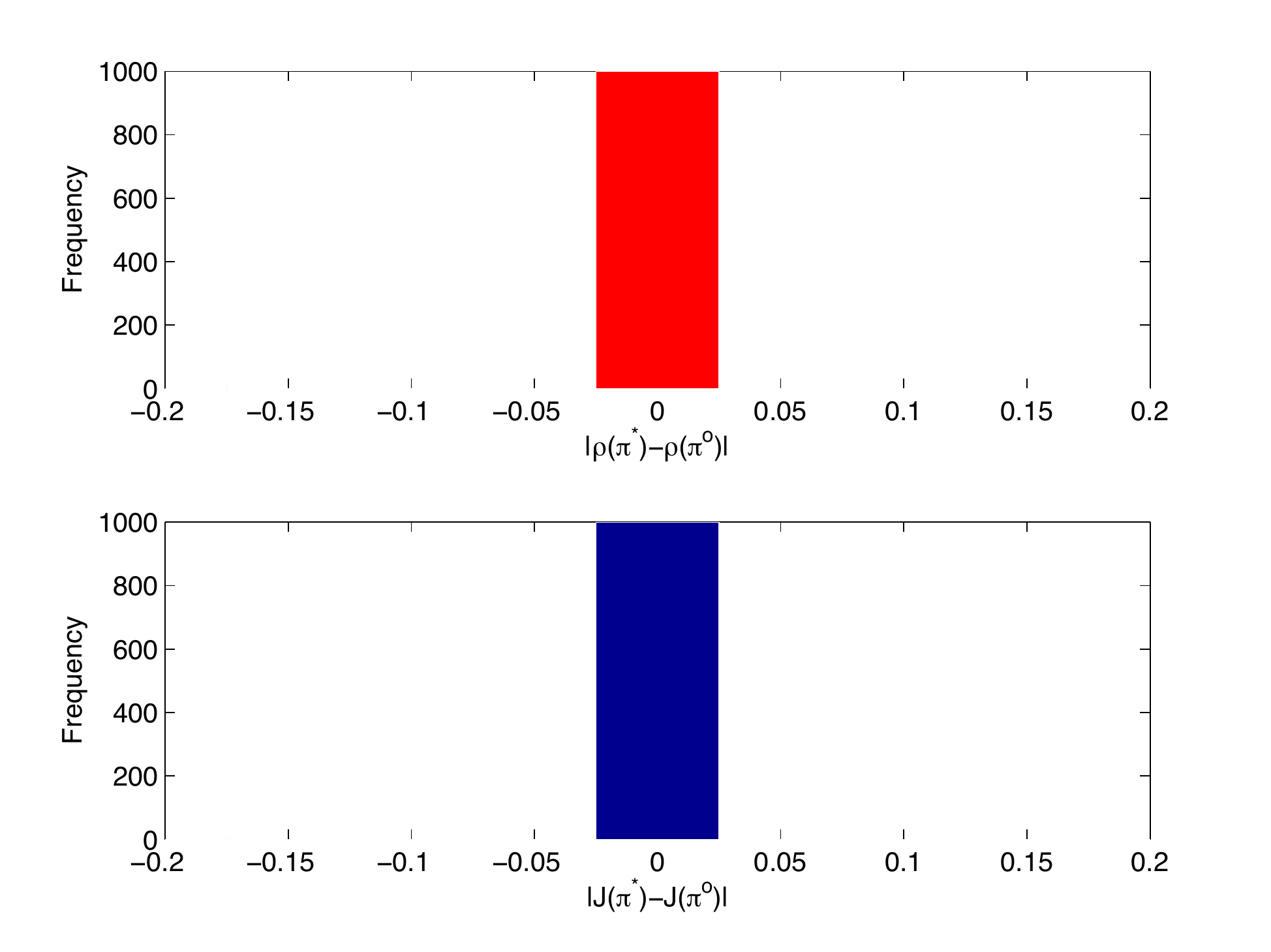} 
      \caption{Histograms of the difference between the average costs corresponding to the optimal and Pareto control policies $\pi^*$ and $\pi^o$ respectively.}
      \label{fig:Pareto}
\end{figure}

We perform 1,000 replications and we observe in Fig. \ref{fig:Pareto} that the absolute difference between $\rho(\pi^0)$ and $\rho(\pi^*)$ is zero. This indicates that $\pi^0$ yields in fact the strong Pareto solution for the one-stage expected costs. Furthermore, Fig. \ref{fig:Pareto} shows that $\min_{\pi\in\Pi}\rho(\pi)=\rho(\pi^0)$, where $\pi^0$ is such that $J^*=J(\pi^0)$. Hence, based on these synthetic data, one can conclude that the optimal control policy of the entire system is the Pareto control policy.

\subsection{Power Management Control of a Hybrid Electric Vehicle: A System with Subsystems of a Principal Group}
\begin{figure}[ht]
  \centering
    \includegraphics[width=3.7 in]{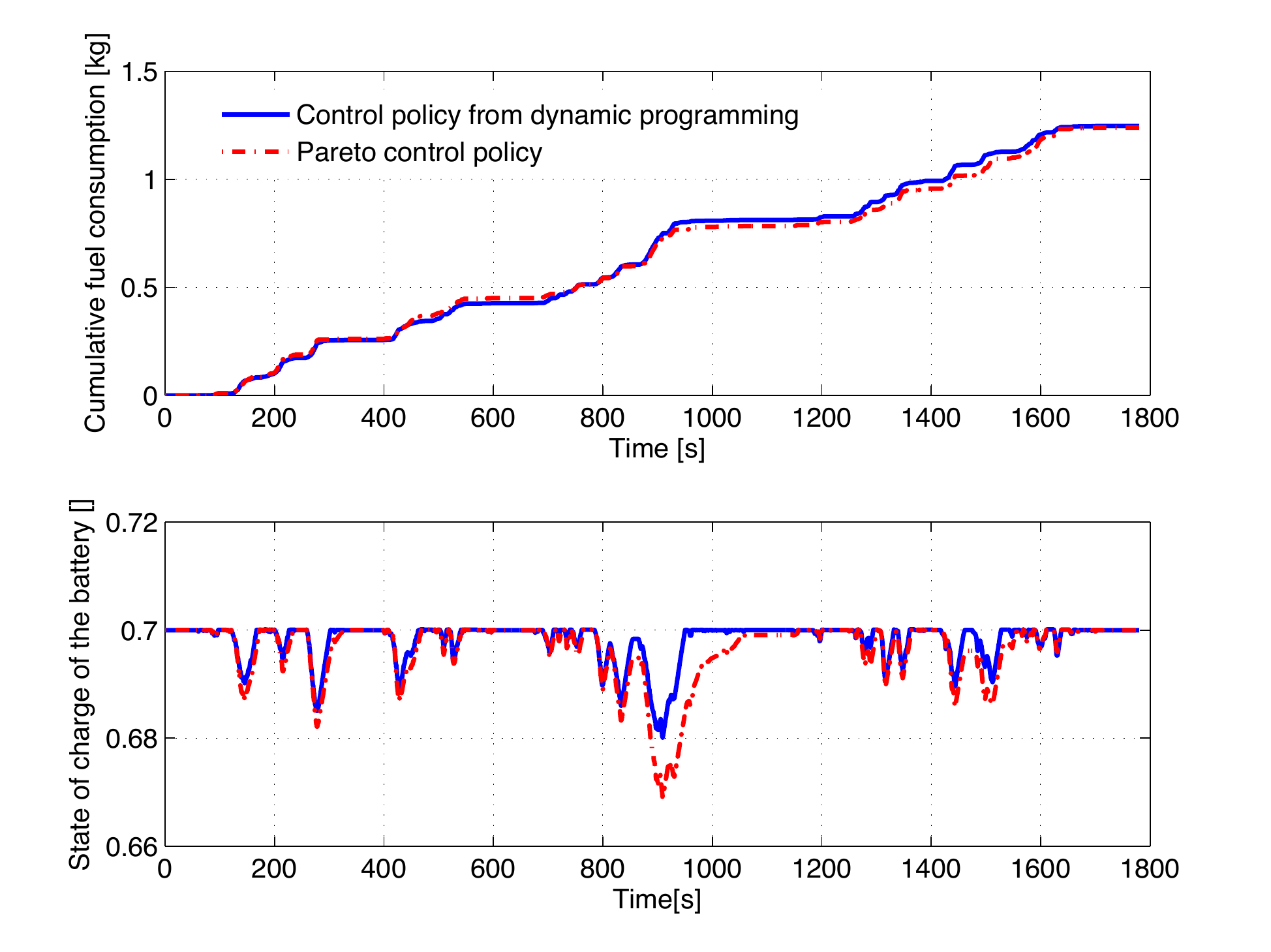} 
      \caption{Cumulative fuel consumption and state of charge of the battery for a parallel hybrid electric vehicle using the control policy derived from dynamic programming and the Pareto control policy over the city-suburban heavy duty vehicle route driving cycle \cite{Malikopoulos2015a}.}
            \label{fig:fuel}
\end{figure}

The results presented here have been used in the problem of optimizing online the power management control in a HEV  \cite{Malikopoulos2015a} consisting of subsystems of a principal group. The Pareto control policy was validated through simulation and it was compared with the control policy derived offline by DP using the long-run expected average cost.  Both control policies achieved the same cumulative fuel consumption as illustrated in Fig. \ref{fig:fuel}, demonstrating that the Pareto control policy is the optimal control policy with respect to the average cost criterion and can be implemented online. This work has been extended \cite{Malikopoulos2014c} by considering the battery in the problem formulation in addition to the engine's and motor's efficiency that can provide insights on how to prioritize these objectives based on consumers' needs and preferences. 

\section{CONCLUDING REMARKS}
In this paper, we established a framework for the analysis and stochastic optimization of complex systems consisting of interactive subsystems. We formulated the stochastic control problem as a multiobjective optimization problem of the one-stage expected costs of the subsystems and developed a duality framework to prove that the Pareto control policy minimizes the long-run expected average cost criterion of the system. We provided a geometric interpretation of the solution and conditions for its existence. The Pareto control policy identifies an equilibrium operating point among the subsystem. If the system operates at this equilibrium, then the long-run expected average cost per unit time is minimized. For practical situations with constraints consistent to those studied here, our results imply that the Pareto control policy may be of value when we seek to derive online the optimal control policy in complex systems. 

One potential extension of this work could be to investigate whether a similar analysis can yield the desired \textit{emergence} in a complex system from a decentralized perspective. Emergence refers to the spontaneous creation of order and functionality from the bottom up. Wherever we see complex systems in the physical world, we see emergent patterns at every level, both in structure and functionality. Emergence occurs without a central planner, from the bottom up, based on the interaction of the individual entities in a system. As a simple example from the natural world of how emergence arises, we can consider the flying patterns created by a flock of birds following three simple rules: 1) stay close but don't bomb into birds around me, 2) fly as fast as birds near me, and 3) move towards the center of the group. The fact that a rule applied locally leads to a macro-level property is what is meant by the term bottom up. Another example of a bottom-up emergent phenomenon is the traffic jam resulting from a specific sequence of vehicle-to-vehicle and vehicle-to-infrastructure interactions. If we could develop the framework to characterize emergence, then we would be able to designate the rules for the interactions of the individual subsystems so that the desired emergent phenomena would occur.

\section{ACKNOWLEDGMENTS}
The author would like to thank Yang Shen for her assistance in running the simulation of the first illustrative example (Section V-B) with varying transition probability and cost matrices.

\bibliographystyle{IEEETran}
\bibliography{AverageCost}

\begin{thebibliography}{10}
\providecommand{\url}[1]{#1}
\csname url@rmstyle\endcsname
\providecommand{\newblock}{\relax}
\providecommand{\bibinfo}[2]{#2}
\providecommand\BIBentrySTDinterwordspacing{\spaceskip=0pt\relax}
\providecommand\BIBentryALTinterwordstretchfactor{4}
\providecommand\BIBentryALTinterwordspacing{\spaceskip=\fontdimen2\font plus
\BIBentryALTinterwordstretchfactor\fontdimen3\font minus
  \fontdimen4\font\relax}
\providecommand\BIBforeignlanguage[2]{{%
\expandafter\ifx\csname l@#1\endcsname\relax
\typeout{** WARNING: IEEEtran.bst: No hyphenation pattern has been}%
\typeout{** loaded for the language `#1'. Using the pattern for}%
\typeout{** the default language instead.}%
\else
\language=\csname l@#1\endcsname
\fi
#2}}

\bibitem{Page20007}
J.~H. Miller and S.~E. Page, \emph{{Complex Adaptive Systems: An Introduction
  to Computational Models of Social Life}}.\hskip 1em plus 0.5em minus
  0.4em\relax Princeton University Press, March 2007.

\bibitem{DOE2011}
F.~Alexander, M.~Anitescu, J.~Bell, D.~Brown, M.~Ferris, M.~Luskin,
  S.~Mehrotra, B.~Moser, A.~Pinar, A.~Tartakovsky, K.~Willcox, S.~Wrigh, and
  V.~Zavala, ``A multifaceted mathematical approach for complex systems,'' DOE
  Workshop on Mathematics for the Analysis, Simulation, and Optimization of
  Complex Systems, Tech. Rep., 2011.

\bibitem{Malikopoulos2014b}
A.~A. Malikopoulos, ``{Supervisory Power Management Control Algorithms for
  Hybrid Electric Vehicles: A Survey},'' \emph{IEEE Transactions on Intelligent
  Transportation Systems}, preprints available online at ieeexplore.ieee.org,
  2014.

\bibitem{Prasanna2013a}
P.~R. Prasanna and A.~Rathore, ``Analysis, design, and experimental results of
  a novel soft-switching snubberless current-fed half-bridge front-end
  converter-based pv inverter,'' \emph{IEEE Transactions on Power Electronics},
  vol.~28, no.~7, pp. 3219--3230, 2013.

\bibitem{Arapostathis1993}
A.~Arapostathis, V.~Borkar, E.~Fernandez-Gaucherand, M.~K. Ghosh, and S.~I.
  Marcus, ``{Discrete-time controlled Markov processes with average cost
  criterion: a survey},'' \emph{SIAM Journal on Control and Optimization},
  vol.~31, no.~2, pp. 282--344, 1993.

\bibitem{Varaiya1978}
P.~Varaiya, ``{Optimal and suboptimal stationary controls for Markov chains},''
  \emph{IEEE Transactions on Automatic Control}, vol. AC-23, no.~3, pp.
  388--394, 1978.

\bibitem{Berstekas2007}
D.~P. Bertsekas and S.~E. Shreve, \emph{{Stochastic Optimal Control: The
  Discrete-Time Case}}, 1st~ed.\hskip 1em plus 0.5em minus 0.4em\relax Athena
  Scientific, February 2007.

\bibitem{Kushner1971}
H.~J. Kushner, \emph{{Introduction to Stochastic Control}}.\hskip 1em plus
  0.5em minus 0.4em\relax Holt, Rinehart and Winston, 1971.

\bibitem{Kumar1986}
P.~R. Kumar and P.~Varaiya, \emph{{Stochastic systems}}.\hskip 1em plus 0.5em
  minus 0.4em\relax Prentice Hall, June 1986.

\bibitem{Howard1960}
R.~A. Howard, \emph{{Dynamic Programming and Markov Processes}}.\hskip 1em plus
  0.5em minus 0.4em\relax The MIT Press, June 1960.

\bibitem{Doob1990}
J.~L. Doob, \emph{{Stochastic Processes}}.\hskip 1em plus 0.5em minus
  0.4em\relax Wiley-Interscience, January 1990.

\bibitem{Malikopoulos2011}
A.~A. Malikopoulos, ``{Equilibrium Control Policies for Markov Chains},'' in
  \emph{50th IEEE Conference on Decision and Control and European Control
  Conference}, Orlando, Florida, {December 12-14} 2011, pp. 7093--7098.

\bibitem{Bellman1957}
R.~Bellman, Ed., \emph{{Dynamic Programming}}.\hskip 1em plus 0.5em minus
  0.4em\relax Princeton University Press, 1957.

\bibitem{White1963}
D.~J. White, ``{Dynamic programming, Markov chains, and the method of
  successive approximations},'' \emph{Journal of Mathematical Analysis and
  Applications}, vol.~6, no.~3, pp. 373--376, 1963.

\bibitem{Bather1973}
J.~Bather, ``{Optimal decision procedures for finite Markov chains. I.
  Examples},'' \emph{Advances in Applied Probability}, vol.~5, no.~2, pp.
  328--339, 1973.

\bibitem{Bather1973a}
------, ``Optimal decision procedures for finite markov chains. part ii:
  Communicating systems,'' \emph{Advances in Applied Probability}, vol.~5,
  no.~3, pp. 521--540, 1973.

\bibitem{Bather1973b}
------, ``Optimal decision procedures for finite markov chains. part iii
  general convex systems,'' \emph{Advances in Applied Probability}, vol.~5,
  no.~3, pp. 541--553, 1973.

\bibitem{Hubner1977}
G.~H\"{u}bner, ``{On the Fixed Points of the Optimal Reward Operator in
  Stochastic Dynamic Programming with Discount Factor Greater than One},''
  \emph{ZAMM - Journal of Applied Mathematics and Mechanics / Zeitschrift
  f\"{u}r Angewandte Mathematik und Mechanik}, vol.~57, no.~8, pp. 477--480,
  1977.

\bibitem{Federgruen1978}
A.~Federgruen, P.~J. Schweitzer, and H.~C. Tijms, ``{Contraction mappings
  underlying undiscounted Markov decision problems},'' \emph{Journal of
  Mathematical Analysis and Applications}, vol.~65, no.~3, pp. 711--730, 1978.

\bibitem{Bertsekas1998}
D.~P. Bertsekas, ``A new value iteration method for the average cost dynamic
  programming problem,'' \emph{SIAM Journal on Control and Optimization},
  vol.~36, no.~2, pp. 742--759, 1998.

\bibitem{Manne1960}
A.~S. Manne, ``{Linear Programming and Sequential Decisions},''
  \emph{Management Science}, vol.~6, no.~3, pp. 259--267, 1960.

\bibitem{Wagner1960}
H.~M. Wagner, ``{On the Optimality of Pure Strategies},'' \emph{Management
  Science}, vol.~6, no.~3, pp. 268--269, 1960.

\bibitem{Derman1962}
C.~Derman, ``{On Sequential Decisions and Markov Chains},'' \emph{Management
  Science}, vol.~9, no.~1, pp. 16--24, 1962.

\bibitem{Derman1965}
C.~Derman and M.~Klein, ``{Some Remarks on Finite Horizon Markovian Decision
  Models},'' \emph{Operations Research}, vol.~13, no.~2, pp. 272--278, 1965.

\bibitem{Hordijk1979}
A.~Hordijk and L.~C.~M. Kallenberg, ``{Linear Programming and Markov Decision
  Chains},'' \emph{Management Science}, vol.~25, no.~4, pp. 352--362, 1979.

\bibitem{Hordijk1984}
------, ``{Constrained Undiscounted Stochastic Dynamic Programming},''
  \emph{Mathematics of Operations Research}, vol.~9, no.~2, pp. 276--289, 1984.

\bibitem{Ross1989a}
K.~W. Ross and R.~Varadarajan, ``{Markov Decision Processes with Sample Path
  Constraints: The Communicating Case},'' \emph{Operations Research}, vol.~37,
  no.~5, pp. 780--790, 1989.

\bibitem{Ross1989b}
K.~W. Ross, ``{Randomized and Past-Dependent Policies for Markov Decision
  Processes with Multiple Constraints},'' \emph{Operations Research}, vol.~37,
  no.~3, pp. 474--477, 1989.

\bibitem{Lasserre1994}
J.~B. Lasserre, ``{Detecting Optimal and Non-Optimal Actions in Average-Cost
  Markov Decision Processes},'' \emph{Journal of Applied Probability}, vol.~31,
  no.~4, pp. 979--990, 1994.

\bibitem{Zadorojniy2006}
A.~Zadorojniy and A.~Shwartz, ``{Robustness of policies in constrained Markov
  decision processes},'' \emph{Automatic Control, IEEE Transactions on},
  vol.~51, no.~4, pp. 635--638, 2006.

\bibitem{Miller1969}
B.~L. Miller and A.~F. {Veinott Jr.}, ``{Discrete Dynamic Programming with a
  Small Interest Rate},'' \emph{The Annals of Mathematical Statistics},
  vol.~40, no.~2, pp. 366--370 CR -- Copyright \&\#169; 1969 Institute of Ma,
  1969.

\bibitem{Chang2007}
H.~S. Chang, ``{A policy improvement method for constrained average Markov
  decision processes},'' \emph{Operations Research Letters}, vol.~35, no.~4,
  pp. 434--438, 2007.

\bibitem{Lamond1989}
B.~F. Lamond and M.~L. Puterman, ``{Generalized Inverses in Discrete Time
  Markov Decision Processes},'' \emph{SIAM Journal on Matrix Analysis and
  Applications}, vol.~10, no.~1, pp. 118--134, 1989.

\bibitem{Ghosh1990}
M.~K. Ghosh, ``Markov decision processes with multiple costs,''
  \emph{Operations Research Letters}, vol.~9, no.~4, pp. 257--260, 1990.

\bibitem{Ren2001}
Z.~Ren and B.~H. Krogh, ``{Adaptive control of Markov chains with average
  cost},'' \emph{Automatic Control, IEEE Transactions on}, vol.~46, no.~4, pp.
  613--617, 2001.

\bibitem{Abounadi2001}
J.~Abounadi, D.~Bertsekas, and V.~S. Borkar, ``Learning algorithms for markov
  decision processes with average cost,'' \emph{SIAM Journal on Control and
  Optimization}, vol.~40, no.~3, pp. 681--698, 2001.

\bibitem{Chang2009}
H.~S. Chang, ``{Decentralized Learning in Finite Markov Chains: Revisited},''
  \emph{Automatic Control, IEEE Transactions on}, vol.~54, no.~7, pp.
  1648--1653, 2009.

\bibitem{Cavazos-Cadena2009}
R.~Cavazos-Cadena, ``{Solutions of the average cost optimality equation for
  finite Markov decision chains: risk-sensitive and risk-neutral criteria},''
  \emph{Mathematical Methods of Operations Research}, vol.~70, no.~3, pp.
  541--566, 2009.

\bibitem{DOE2006}
B.~Hendrickson and M.~Wright, ``Mathematical research challenges in
  optimization of complex systems,'' Sandia National Laboratories and Courant
  Institute of Mathematical Sciences, Tech. Rep., 2006.

\bibitem{Grimmett2001}
G.~R. Grimmett and D.~R. Stirzaker, \emph{{Probability and Random Processes}},
  3rd~ed.\hskip 1em plus 0.5em minus 0.4em\relax Oxford University Press,
  August 2001.

\bibitem{Ross1995}
S.~M. Ross, \emph{{Stochastic Processes}}, 2nd~ed.\hskip 1em plus 0.5em minus
  0.4em\relax Wiley, January 1995.

\bibitem{Ehjrgott2005}
M.~Ehrgott, \emph{Multicriteria Optimization}.\hskip 1em plus 0.5em minus
  0.4em\relax Springer, 2nd edition, 2005.

\bibitem{Bertsekas2003}
D.~P. Bertsekas, A.~Nedic, and A.~E. Ozdaglar, \emph{{Convex Analysis and
  Optimization}}.\hskip 1em plus 0.5em minus 0.4em\relax Athena Scientific,
  April 2003.

\bibitem{Searle2006}
S.~Searle, G.~Casella, and C.~McCulloch, \emph{Variance Components}, ser. Wiley
  Series in Probability And Statistics.\hskip 1em plus 0.5em minus 0.4em\relax
  Wiley, 2006.

\bibitem{Horn1994}
R.~Horn and C.~Johnson, \emph{Topics in Matrix Analysis}, ser. Topics in Matrix
  Analysis.\hskip 1em plus 0.5em minus 0.4em\relax Cambridge University Press,
  1994.

\bibitem{Malikopoulos2015}
A.~A. Malikopoulos, V.~Maroulas, and J.~Xiong, ``A multiobjective optimization
  framework for stochastic control of complex systems,'' in \emph{Proceedings
  of the 2015 American Control Conference}, 2015, pp. 4263--4268.

\bibitem{Malikopoulos2015a}
A.~A. Malikopoulos, ``A multiobjective optimization framework for online
  stochastic optimal control in hybrid electric vehicles,'' \emph{IEEE
  Transactions on Control Systems Technology}, 2015 (forthcoming).

\bibitem{Malikopoulos2014c}
M.~Shaltout, A.~A. Malikopoulos, S.~Pannala, and D.~Chen, ``A consumer-oriented
  control framework for performance analysis in hybrid electric vehicles,''
  \emph{IEEE Transactions on Control Systems Technology}, vol.~23, no.~4, pp.
  1451--1464, 2015.

\end{thebibliography}

\begin{IEEEbiography}[{\includegraphics[width=1.1in,height=1.25in,clip,keepaspectratio]{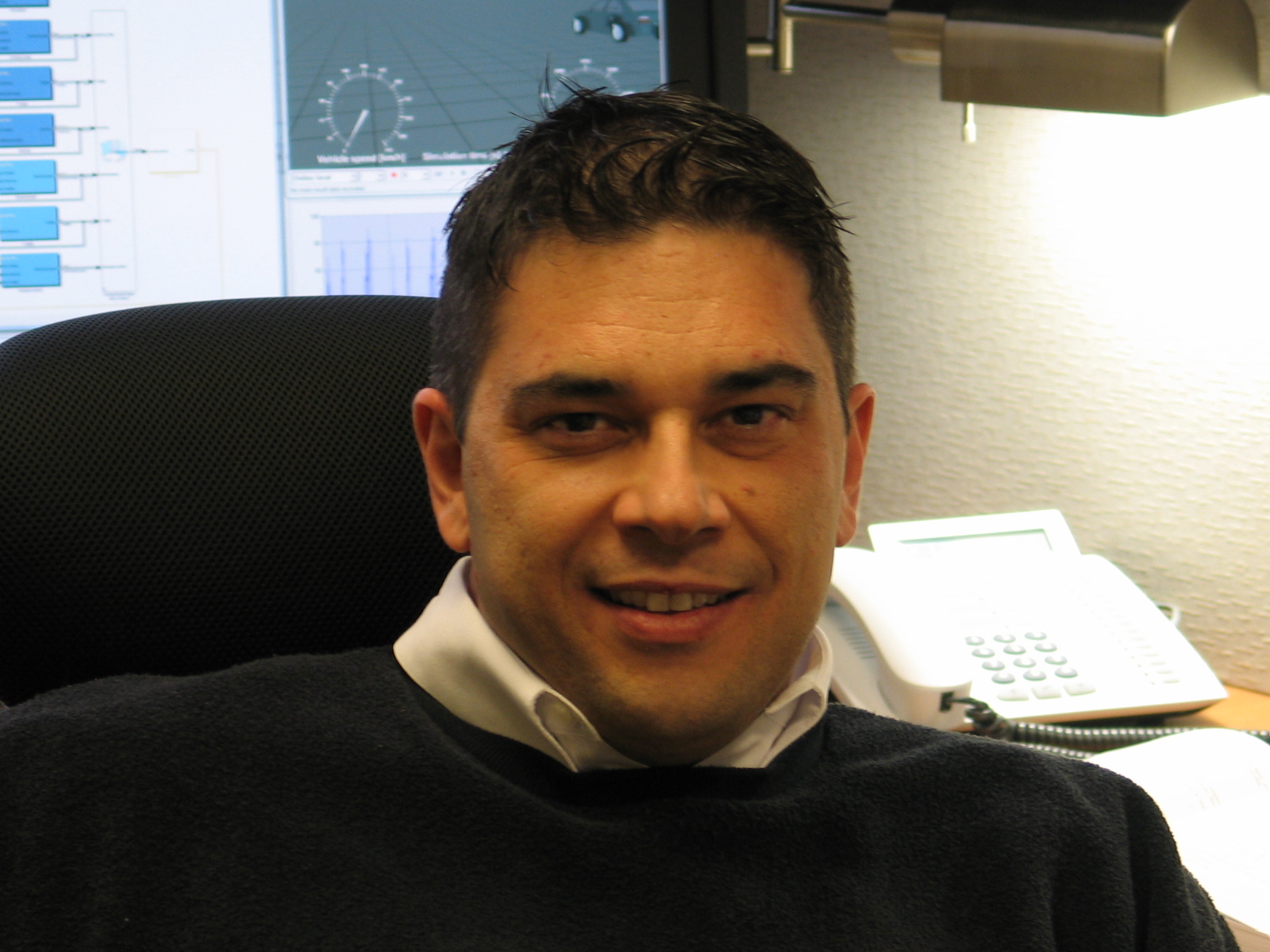}}]{Andreas A. Malikopoulos}
 (M2006) received a Diploma in Mechanical Engineering from the National Technical University of Athens, Greece, in 2000. He received M.S. and Ph.D. degrees from the Department of Mechanical Engineering at the University of Michigan, Ann Arbor, Michigan, USA, in 2004 and 2008, respectively.

He is currently the Deputy Director of the Urban Dynamics Institute and an Alvin M. Weinberg Fellow with the Energy \& Transportation Science Division at Oak Ridge National Laboratory (ORNL). Before joining ORNL he was a Senior Researcher with General Motors Global Research \& Development, conducting research in the area of stochastic optimization and control of advanced propulsion systems. His research spans several fields, including analysis, optimization, and control of complex systems; decentralized systems; and stochastic scheduling and resource allocation problems. The emphasis is on applications related to energy, transportation, and operations research.
\end{IEEEbiography}

\end{document}